\documentclass{article}

\usepackage{amsmath,amsthm,amsfonts, mathrsfs, amssymb,latexsym,epsfig}
\usepackage{url}

\theoremstyle{plain}
\newtheorem{theorem}{Theorem}
\newtheorem{lemma}{Lemma}
\newtheorem{corollary}{Corollary}
\begin{document}
\title{Irrationality proofs \`a la Hermite}
\author{Li Zhou}
\date{}
\maketitle

\section{Introduction.}

In \cite{n} Niven used the integral $$H_n=\int_0^{\pi}\frac{x^n(\pi -x)^n}{n!}\sin x\, dx$$ to give a well-known proof of the irrationality of $\pi$. Recently Zhou and Markov \cite{zm} used a recurrence relation satisfied by $H_n$ to present an alternative proof which may be more direct than Niven's.

Niven did not cite any reference in \cite{n} and thus the origin of $H_n$ seems rather mysterious and ingenious. However if we heed Abel's advice to ``read the masters'', we find that $H_n$ emerged much more naturally from the great works of Lambert \cite{l} and Hermite \cite{h}. In fact, we find that Hermite had already used $H_n/2^{n+1}$ to give a simple proof of the irrationality of $\pi^2$. The re-examination of Hermite's works also leads us to a short new proof of the irrationality of $r\tan r$ for $r^2\in\mathbb{Q}\setminus\{0\}$, and a generalisation to the irrationality of certain ratios of Bessel functions. 

\section{The origin of $H_n$.} 

It is well known that in 1761 Lambert conceived the first proof of the irrationality of $\tan r$ for nonzero rational $r$, and as a corollary, the irrationality of $\pi$. Lambert started with 
$$\tan r=\frac{\sin r}{\cos r}=\frac{r-\frac{r^3}{3!}+\frac{r^5}{5!}-\cdots}{1-\frac{r^2}{2!}+\frac{r^4}{4!}-\cdots}$$ and used the Euclidean algorithm to construct the continued fraction $$\tan r=\frac{r}{1-\frac{r^2}{3-\frac{r^2}{5-\frac{r^2}{\ddots}}}}$$ and its related ``remainders'' $R_1, R_2, R_3, \ldots $, where
\begin{eqnarray*}
R_1 &=& \sin r-r\cos r = \frac{r^3}{3}-\frac{r^5}{2\cdot 3\cdot 5}+\cdots,\\
R_2 &=& (3-r^2)\sin r -3r\cos r= \frac{r^5}{3\cdot 5}-\frac{r^7}{2\cdot 3\cdot 5\cdot 7}+\cdots, \\ 
R_3 &=& (15-6r^2)\sin r-(15r-r^3)\cos r= \frac{r^7}{3\cdot 5\cdot 7}-\frac{r^9}{2\cdot 3\cdot 5\cdot 7\cdot 9}+\cdots, \\
&\ldots & .\end{eqnarray*}
Intuitively, these remainders measure how closely $\tan r$ is approximated by the truncations of its continued fraction. For example, $R_2$ measures the closeness of $$\frac{\sin r}{\cos r}\approx\frac{r}{1-\frac{r^2}{3}}=\frac{3r}{3-r^2}.$$ After much labor Lambert accomplished his feat by studying the recurrence and convergence properties of these remainders.

Less well known is that in 1873 Hermite was attracted to these remainders and noticed that they satisfy two differential relations 
\begin{equation}\label{e01}\frac{dR_n}{dr}=rR_{n-1} ~\textrm{and}~ \frac{d}{dr}\left(\frac{R_{n-1}}{r^{2n-1}}\right)=-\frac{R_n}{r^{2n}}.\end{equation} 
From (\ref{e01}) he further derived Lambert's recurrence relation  $R_{n}=(2n-1)R_{n-1}-{r}^2 R_{n-2}$, a differential equation 
\begin{equation}\label{e02}\frac{d^2R_n}{dr^2}-\frac{2n}{r}\frac{dR_n}{dr}+R_n=0, \end{equation} and an integral representation \begin{equation}\label{e03}R_n(r)=\frac{r^{2n+1}}{2^n n!}\int_0^1 (1-z^2)^n\cos (rz) \,dz.\end{equation} Using this integral, with $r=\pi/2$, Hermite then gave a simple proof of the irrationality of $\pi^2$.

Now substituting $z=1-2x/\pi $ and $r=\pi/2$ into (\ref{e03}), we have
$$R_n\left(\frac{\pi}{2}\right)=\frac{1}{2^n n!}\int_0^{\pi/2} x^n(\pi-x)^n\sin x \,dx,$$ which is of course $H_n/2^{n+1}$ by symmetry. Hence Niven's widely-known simple proof in \cite{n} is neither that distant from Lambert's original idea nor that different from Hermite's \emph{already-simple} proof of a stronger result. Pedagogically, it is perhaps more regrettable that Hermite's proof has rarely been acknowledged since the publication of \cite{n}, because Hermite presented his ingenious ideas without covering up their origins and motivations. 

\section{Two possible paths to Hermite's integral.}

From (\ref{e01}) it is a simple exercise to derive Lambert's recurrence relation and (\ref{e02}), but not so easy to obtain (\ref{e03}). Hermite omitted the derivation of (\ref{e03}), either because he considered it routine or because he was aware of something else which we shall mention later. For now we present two natural paths from (\ref{e01}) or (\ref{e02}) to (\ref{e03}) since they may be motivating and interesting to teachers and students of calculus and differential equations.

For the first approach we rewrite the first relation in (\ref{e01}) as $ R_n = \int_0^r {tR_{n-1}(t)} dt $ and use it to integrate inductively. To be precise, we start with any integrable function $f(x)$ and $$R_0=\int_0^r f(t) dt = r\int_0^1 f(rz) dz.$$ Suppose that we have derived for some $n\ge 1$, 
\begin{eqnarray*}R_{n-1} &=& \int_0^{t_{n}=r}t_{n-1}\int_0^{t_{n-1}} \cdots t_{1}\int_0^{t_1}f(t_0) dt_0  \cdots dt_{n-2} dt_{n-1} \\
&=& \frac{r^{2n-1}}{2^{n-1} (n-1)!}\int_0^1 (1-z^2)^{n-1}f(rz) dz.\end{eqnarray*} Then applying this to $g(x)=x\int_0^x f(t) dt$ we get
\begin{eqnarray*}
R_n &=& \int_0^{t_{n+1}=r} t_{n}\int_0^{t_{n}} \cdots t_{1}\int_0^{t_{1}}f(t_{0}) dt_{0} \cdots dt_{n-1} dt_{n}\\
&=& \int_0^{t_{n+1}=r} t_{n}\int_0^{t_{n}} \cdots t_{2}\int_0^{t_{2}}g(t_1) dt_{1}  \cdots dt_{n-1} dt_{n}\\
&=& \frac{r^{2n-1}}{2^{n-1} (n-1)!}\int_0^1 (1-z^2)^{n-1}g(rz) dz\\
&=& \frac{r^{2n-1}}{2^{n-1}(n-1)!}\int_0^1 (1-z^2)^{n-1} rz\int_0^{rz}f (t) dt dz\\
&=& \frac{r^{2n}}{2^{n-1}(n-1)!}\left(\left[\frac{(1-z^2)^n}{-2n}\int_0^{rz} f(t) dt\right]_{z=0}^{z=1}+\int_0^1 \frac{(1-z^2)^n}{2n}r f(rz) dz\right)\\
&=& \frac{r^{2n+1}}{2^n n!}\int_0^1 (1-z^2)^n f(rz) dz.\end{eqnarray*} Letting $f(x)=\cos x$ completes the first derivation. 

For the second approach we consider $r\ge 0$ and substitute $r=\sqrt{t}$ and $y_n(t)=R_n(\sqrt{t})$ into (\ref{e02}) to obtain
$$ 4t\frac{d^2y_n}{dt^2}-(4n-2)\frac{dy_n}{dt}+y_n = 0; ~y_n(0)=0.$$
Taking Laplace transform $\mathscr{L}$ we get 
$$4s^2\frac{dY_n}{ds}+[(4n+6)s-1]Y_n = 0, $$ where $Y_n(s)=\mathscr{L}\{y_n(t)\}(s)$. Hence $$Y_n(s)=C_n s^{-n-\frac{3}{2}}e^{-\frac{1}{4s}} = \frac{C_n}{n!}\left(\frac{n!}{s^{n+1}}\cdot \frac{e^{-\frac{1}{4s}}}{\sqrt{s}}\right),$$ where $C_n$ is a constant. From a respectable table of Laplace transforms \cite{p} we find $$\mathscr{L}\left\{\frac{\cos\sqrt{t}}{\sqrt{t}}\right\}=\sqrt{\frac{\pi}{s}}e^{-\frac{1}{4s}}.$$ Therefore by the convolution theorem,
$$y_n(t)=\frac{C_n}{n!\sqrt{\pi}}\left(t^n\ast\frac{\cos\sqrt{t}}{\sqrt{t}}\right)=\frac{C_n}{n!\sqrt{\pi}}\int_0^t (t-v)^n\frac{\cos\sqrt{v}}{\sqrt{v}} dv .$$ The substitutions $v=x^2$ and $t=r^2$ yield 
$$R_n(r)=\frac{2C_n}{n!\sqrt{\pi}}\int_0^r(r^2-x^2)^n\cos x dx=\frac{2C_n r^{2n+1}}{n!\sqrt{\pi}}\int_0^1 (1-z^2)^n\cos(rz) dz .$$ Noticing that 
$$ \frac{R_n}{r^{2n+1}}=\frac{1}{3\cdot 5\cdot 7\cdots (2n+1)}$$ at $r=0$, we get
$$\frac{2C_n}{n!\sqrt{\pi}}\int_0^1(1-z^2)^n dz=\frac{1}{3\cdot 5\cdot 7\cdots (2n+1)}$$ from which it is then easy to figure out that $C_n=\sqrt{\pi}/2^{n+1}$.

\section{A new proof of irrationality.}

Another reward of reading the masters is that we can also use Hermite's integral to give a new, simple, and self-contained proof of the irrationality of $\tan r$ for nonzero rational $r$. In fact we can do slightly better. 

\begin{theorem}\label{th1} If $r^2\in \mathbb{Q}\setminus\{0\}$ then $r\tan r$ is irrational.
\end{theorem}

\begin{proof} The irrationality of $\pi^2$ will be a byproduct of this proof, so we start by assuming that $r^2\in \mathbb{Q}\setminus\{0\}$ and $\cos r\ne 0$. Write $r^2=a/b$ with $a,b\in\mathbb{Z}$ and assume that $r\tan r=p/q$ with $p,q\in \mathbb{Z}$. For $n\ge 0$, let $$f_n(x)=\frac{(r^2-x^2)^n}{2^n n!} ~\textrm{and}~ R_n=\int_0^r f_n(x)\cos x \, dx.$$ Then $b^{\left\lceil n/2\right\rceil}R_n \to 0$ as $n\to \infty$, $R_0=\sin r$, and $R_1=\sin r-r\cos r$. For $n\ge 2$, it is easy to verify that $f_n''(x)=-(2n-1)f_{n-1}(x)+r^2f_{n-2}(x)$. Integrating by parts twice we then have 
\begin{equation}\label{e04}
R_{n}=-\int_0^r f_n''(x)\cos x\,dx=(2n-1)R_{n-1}-{r}^2 R_{n-2}.\end{equation}
Induction on $n$ in (\ref{e04}) shows that for $n\ge 0$, $R_n=u_n\sin r+v_n\cos r$ where $u_n$ and $rv_n$ are polynomials in $r^2$ with integer coefficients and degrees at most $\left\lceil  n/2\right\rceil$. %Thus \begin{equation}\label{e05}\frac{qrb^{\left\lceil n/2\right\rceil}R_n }{\cos r}=b^{\left\lceil n/2\right\rceil}(u_np+rv_nq)\end{equation} which is an integer for all $n\ge 0$.%
Moreover, if two consecutive terms of the sequence $\langle  R_n\rangle$ are $0$, then (\ref{e04}) forces all terms of $\langle  R_n\rangle$ to be $0$, contradicting the fact that $R_0-R_1=r\cos r\ne 0$. Hence $\langle R_n \rangle$ has infinitely many nonzero terms. %Moreover, $$|R_n|\le \frac{|r^{2n+1}|}{2^n n!}, ~\textrm{so}~ |b^{\left\lceil n/2\right\rceil}R_n|\le |b^nR_n|\le \frac{|r(a/2)^n|}{n!}\to 0 ~\textrm{as}~ n\to \infty.$$ Therefore,% 
Therefore we can pick a large enough $n$ such that $qrb^{\lceil n/2\rceil}R_n/\cos r=b^{\lceil n/2\rceil}(u_np+rv_nq)$ is a nonzero integer in $(-1,1)$, a contradiction. 

Since $\pi\tan\pi=0\in\mathbb{Q}$ and $\cos \pi=-1\ne 0$, $\pi^2\notin\mathbb{Q}$.  Thus the condition $r^2\in \mathbb{Q}\setminus\{0\}$ automatically implies that $\cos r\ne 0$. Therefore we have proved that $r\tan r\notin\mathbb{Q}$ whenever $r^2\in \mathbb{Q}\setminus\{0\}$.
\end{proof}

This proof fully showcases the advantage of Hermite's integral approach: $R_n$ is easy to define as an integral in a self-contained manner; the limiting property of $R_n$ is immediate; and the recurrence relation satisfied by $R_n$ is a simple consequence of integration by parts. It is thus not surprising that the popular modern proofs of the irrationality of $\pi$ and $\pi^2$ are either slight variations or rediscoveries of Hermite's original one (for example, see \cite{n}, \cite{zm}, and \cite[pp.\ 117--118]{bm}). However, it has not been noticed until in \cite{zm} that the recurrence relation has an added bonus in establishing the existence of a nonzero subsequence, since in the special case of $r=\pi/2$ the integral $R_n(\pi/2)$ is manifestly positive, so there is no such need. 

The observant reader may also notice that our proof can easily accommodate the case of $r^2<0$, since $f_n(x)\cos x$ is an entire function and thus its integral from $0$ to $r$ is path-independent.  Therefore Theorem 1 includes the implicit statement that $r\tanh r$ is irrational whenever $r^2\in \mathbb{Q}\setminus\{0\}$. An immediate corollary of this is that $e^r$ is irrational for nonzero rational $r$. 

\section{Generalising to Bessel functions.} 

So what attracted Hermite to Lambert's remainders and how did he  ``notice'' the differential relations in (\ref{e01})? The answer may lie in the fact that Hermite referred (\ref{e02}) fleetingly as a Bessel differential equation. Indeed if we change variables by $R_n(r)=r^{n+1/2}w(r)$ then (\ref{e02}) becomes
\begin{equation}\label{e05}r^2\frac{d^2 w}{d r^2}+r\frac{d w}{d r}+\left[r^2-\left(n+\frac{1}{2}\right)^2\right]w=0,\end{equation}
which is the more familiar form of the Bessel equation of order $n+1/2$. As a consequence of this realisation, our second derivation of (\ref{e03}) offers a method of solving the Bessel equation of order $\nu$ not seen in typical textbooks of differential equations. A solution to (\ref{e05}) is $J_{n+1/2}(r)$ where \begin{equation}\label{j01} J_{\nu}(r)=\sum_{k=0}^{\infty} \frac{(-1)^k}{k!\Gamma(\nu+k+1)}\left(\frac{r}{2}\right)^{\nu+2k}\end{equation} is the Bessel function of the first kind of order $\nu$ \cite[Chapter XVII]{ww}. Now, by comparing the expansion of $R_n$ and using the fact that $\Gamma(1/2)=\sqrt{\pi}$, we see that $$R_n(r)=\sqrt{\frac{\pi}{2}}r^{n+\frac{1}{2}}J_{n+\frac{1}{2}}(r).$$ 

This connection with Bessel functions leads us naturally to a generalisation of Theorem \ref{th1}. We start by recalling from \cite{ww} the well-known relations 
\begin{equation}\label{e06}
rJ_{\nu+1}=2\nu J_{\nu}-rJ_{\nu-1}, \end{equation}
\begin{equation}\label{j02} \frac{d}{rdr}\left(r^{-\nu}J_{\nu}(r)\right)=-r^{-(\nu+1)}J_{\nu+1}(r), \end{equation} and Poisson's integral representation 
\begin{equation}\label{e07}J_{\nu}(r)=\frac{r^{\nu}}{2^{\nu}\sqrt{\pi}\Gamma(\nu+\frac{1}{2})}\int_0^{\pi} \cos (r\cos\theta)\sin^{2\nu}\theta \,d\theta; ~\textrm{Re} (\nu)>-\frac{1}{2}. \end{equation}
Notice that (\ref{e03}) follows immediately by letting $z=\cos \theta$ in (\ref{e07}), which is likely the reason why Hermite omitted the derivation of (\ref{e03}). To make our proof of the generalisation cleaner we present a lemma first.

\begin{lemma}\label{l1} For fixed $r\ne 0$ and $\nu$ the sequence $\left\langle J_{\nu+n}(r)\right\rangle_{n\in\mathbb{Z}}$ cannot contain two consecutive zeros.
\end{lemma}

\begin{proof} Suppose that $J_{\nu+m}(r)=J_{\nu+m+1}(r)=0$ for some $m\in\mathbb{Z}$. Then the recurrence relation (\ref{e06}) forces $J_{\nu+n}(r)=0$ for all $n\in\mathbb{Z}$. Using (\ref{j02}) inductively we then have $y_{\nu}^{(n)}(r)=0$ for all $n\ge 0$, where $y_{\nu}(r)=r^{-\nu}J_{\nu}(r)$. Thus $y_{\nu}\equiv 0$, a contradiction. \end{proof} 
 
\begin{theorem}\label{th2} If $s\in\mathbb{Q}$, $r^2\in\mathbb{Q}\setminus\{0\}$, and $J_s(r)\ne 0$, then $rJ_{s+1}(r)/J_s(r)$ is irrational.
\end{theorem}

\begin{proof} Suppose that $r$ and $s$ satisfy the hypothesis. Write $r^2=a/b$ and $s=c/d$  with $a, b, c, d\in\mathbb{Z}$. Assume that $rJ_{s+1}/J_s=p/q$ with $p,q\in \mathbb{Z}$. Induction on $n$ in (\ref{e06}) shows that for $n\ge 0$, $$r^nJ_{n+s+1}=u_nJ_{s+1}+v_nJ_{s}$$ where $u_n, rv_n \in \mathbb{Z} [s, r^2]$ with degrees at most $n$ in $s$ and at most $\left\lceil  n/2\right\rceil$ in $r^2$. Thus \begin{equation}\label{e08}\frac{qb^{\left\lceil n/2\right\rceil}d^nr^{n+1}J_{n+s+1}}{J_s}=b^{\left\lceil n/2\right\rceil}d^n\left(u_np+rv_nq\right)\end{equation} which is an integer for all $n\ge 0$. Also by Lemma \ref{l1}, the sequence $\langle J_{n+s+1} \rangle$ has infinitely many nonzero terms. Moreover, for all large enough $n$, $n+s+1>-1/2$, so (\ref{e07}) yields  $$|J_{n+s+1}|\le \frac{|r^{n+s+1}|\sqrt{\pi}}{2^{n+s+1}\Gamma(n+s+\frac{3}{2})}.$$ Hence we can pick a large enough $n$ such that the expression in (\ref{e08}) is a nonzero integer in $(-1,1)$, a contradiction. \end{proof}

\begin{corollary}\label{c1} If $s\in\mathbb{Q}$ and $r^2\in\mathbb{Q}\setminus\{0\}$ then $J_s(r)\ne 0$.\end{corollary}

\begin{proof} Suppose that $r$ and $s$ satisfy the hypothesis. If $J_s(r)=0$ then $J_{s-1}(r)\ne 0$ by Lemma \ref{l1}, thus $rJ_s(r)/J_{s-1}(r)=0\in\mathbb{Q}$, contradicting Theorem \ref{th2}. \end{proof}

\begin{corollary}\label{c2} If $s\in\mathbb{Q}$ and $r^2\in\mathbb{Q}\setminus\{0\}$ then $rJ_{s+1}(r)/J_s(r)$ is irrational.\end{corollary}

\begin{proof} This follows immediately from Theorem \ref{th2} and Corollary \ref{c1}. \end{proof}

Imitating the analogy between $\tan r$ and $\tanh r$, we can similarly replace $J_{\nu}$ above by $I_{\nu}$, where $I_{\nu}(r)$ is the modified Bessel function of the first kind, defined by
$$I_{\nu}(r)=i^{-\nu}J_{\nu}(ir)=\sum_{k=0}^{\infty} \frac{1}{k!\Gamma(\nu+k+1)}\left(\frac{r}{2}\right)^{\nu+2k}.$$

\paragraph{Acknowledgment.} I am very grateful to my good friend Dr. Lubomir Markov whose nice talk at the 2009 FL-MAA Conference sparked my  interest in the proofs of irrationality. Since then he has generously provided valuable references, stimulating discussions, and constant encouragements to which I express my sincere thanks! I would also like to acknowledge the numerous helpful suggestions of the referee for improving the presentation of this note.

\bigskip

\noindent\textit{Department of Mathematics, Polk State College,
Winter Haven, FL 33881, USA\\
lzhou@polk.edu}\end{document}